\documentclass[twocolumn]{autart}    

\usepackage{mathrsfs}

\usepackage{cite}
\usepackage{amsmath,amssymb,amsfonts}
\usepackage{algorithmic}
\usepackage{graphicx}
\usepackage{textcomp}
\usepackage{mathrsfs}
\usepackage{graphicx}          
\usepackage{subfig}
\usepackage{amstext}
\usepackage{cite}

\newtheorem{lemma}{Lemma}
\newtheorem{definition}{Definition}
\newtheorem{theorem}{Theorem}
\newtheorem{remark}{Remark}

\newtheorem{example}{Example}

\newtheorem{corollary}{Corollary}                               

\begin{document}

\begin{frontmatter}

\title{Discernibility of topological variations for networked LTI systems based on observed output trajectories} 

\thanks[footnoteinfo]{This work is supported by the National Natural Science Foundation of China under Grants 12172020, 11932003, T2121002, in part by the Beijing Natural Science Foundation
under Grant 1222010, and in part by the Hong Kong Research Grants Council under the GRF Grant CityU 11206320. The material in this paper was not presented at any conference.}

\author[YUQING]{Yuqing Hao}\ead{haoyq@buaa.edu.cn},    
\author[YUQING]{Qingyun Wang}\ead{nmqingyun@163.com},
\author[Duan]{Zhisheng Duan}\ead{duanzs@pku.edu.cn},             
\author[CHEN]{Guanrong Chen}\ead{eegchen@cityu.edu.hk}  

\address[YUQING]{Department of Dynamics and Control, Beihang University, Beijing, 100191, China}  
\address[Duan]{State Key Laboratory for Turbulence and Complex Systems, Department of Mechanics and Engineering Science, College of Engineering, Peking University, Beijing 100871, China}             
\address[CHEN]{Department of Electronic Engineering, City University of Hong Kong, Hong Kong, China}        

\begin{keyword}                           
Networked systems, topological variation, discernibility, output trajectory.               
\end{keyword}                             

\begin{abstract}                          
In this paper, the possibility of detecting topological variations by observing output trajectories from networked linear time-invariant systems is investigated, where the network topology can be general{,} but the nodes have identical higher-dimensional dynamics. A necessary and sufficient condition on the discernibility of topological variations is derived, in terms of the eigenspaces of the original and the modified network configuarations. By taking the specific network structures into consideration, some lower-dimensional conditions are derived, which reveal how the network topologies, sensor locations, node-system dynamics and output, as well as inner interactions altogether affect the discernibility. Furthermore, the output discernibility of topological changes for networked multi-agent systems is revisited, showing that some criterion reported in the literature does not hold. Consequently, a modified necessary and sufficient condition is established. The effectiveness of the results is demonstrated through several examples.
\end{abstract}

\end{frontmatter}

\section{Introduction}
In the last two decades, the study of networked systems has gained enormous popularity from the communities in engineering, information technology, mathematics, sociology, biology, and physics. Examples of networked systems in practical applications include wireless communication networks \cite{wood}, networked robotics \cite{bullo}, global transportation networks \cite{banavar}, power generation and distribution networks \cite{wenguanghui}, and so on.

It is well known that topological structure has significant effects on the performance of the underlying networked systems. In practical networks, the failures of network components or denial-of-service attacks may result in topological variations \cite{albert,zuozhiqiang}, which can affect the network performance \cite{wen_switching,wang_controllability_2015,hao_tac,hao_role}, sometimes even have disastrous impacts on the secure and reliable operation. One example is the catastrophic power outage in southern Italy in 2003, which was reportedly caused by failures of some high-voltage transmission lines \cite{buldyrev}. Therefore, there is a need of detecting topological variations in time to protect the networked systems.

Detecting topological variations of networked systems has received compelling attention, with many efficient approaches proposed \cite{rahimian_detection,dhal_detecting,torres,davoodi,costanzo,pandey}. A method to detect and isolate link failures was developed in \cite{rahimian_detection}, which is based on the observed jumps in the derivatives of the output responses of a subset of nodes in networked linear time-invariant (LTI) systems. In \cite{dhal_detecting} and \cite{torres}, link failures in a network synchronization process were detected from noisy local measurements by using a maximum a posteriori probability detection technique. A distributed fault detection and isolation filter was designed in \cite{davoodi} for a network of heterogeneous multi-agent systems. A structural anomaly detection algorithm was proposed in \cite{costanzo}, which is capable of detecting topological changes in dynamical networks. In \cite{pandey}, a diffusion protocol for networked multi-agent systems was constructed and applied to link failure detection. Note that most efforts are devoted to developing detection algorithms. However, the possibility of detecting a topological change by observing the network behavior is even more important but challeging \cite{indiscernible}. 

Recently, the discernibility of topological variations has become a focal topic for investigation.
{The detectability of single link failures in a multi-agent network was investigated in \cite{rahimian_characterization},} which was generalized to multi-link failures in \cite{rahimian_detectability}. {In \cite{battistelli_CDC}, some conditions for detecting node or link disconnections of integrator networks were established. The effects of an edge or a node disconnection on a multi-agent consensus network were investigated in \cite{valcher}. The results were further extended to the detection and identification of edge disconnections in \cite{arxiv}. Note that most if not all existing results on the discernibility of topological variations are derived under the assumption that all nodes are one-dimensional \cite{battistelli_CDC,valcher,arxiv}. However, in real-world networks, nodes typically have higher-dimensional
states, which are coupled via multi-dimensional communication
channels \cite{chenguanrong}. In such
situations, the discernibility of topological changes becomes much more complicated and
challenging.

To date, there has been little work on the discernibility of topological changes for networked higher-dimensional systems. The discernibility of topological changes for networks of linear dynamical systems was studied in \cite{detecting_topology}. A necessary and sufficient condition on the discernibility of topological changes for a network of differential-algebraic systems was established in \cite{indiscernible}, with the indiscernible initial states characterized.} Very recently, some lower-dimensional conditions on the discernibility of topological variations for networked LTI systems were established in \cite{haoyuqing}. Note that all the above-mentioned works focus on detecting topological changes from the observation of the whole network state. However, complete observation of the full state is always unrealistic. In many practical situations, only partial information about the node-system state is accessible, and only a subset of nodes are available for measurement. Hence, detecting topological variations by observing output trajectories of networks has broader applicability in practice. The output discernibility of topological variations for networks of linear dynamical systems was investigated in \cite{detecting_topology}, which requires the network topology to be undirected. In \cite{zhangyuan}, the detectability and isolability of topology failures for a network based on the observed output behavior were studied, where the edge weights can be unknown.

{In this paper, the possibility of detecting topological variations by observing output trajectories of networked LTI systems is investigated. The contribution of this paper is four-fold. First, the network topology can be general, directed and weighted. The node-systems have identical higher-dimensional linear dynamics. Differing from \cite{detecting_topology,indiscernible} and \cite{valcher}, this paper allows directed network topologies. Second, a necessary and sufficient condition on the discernibility of topological variations is derived, in terms of the eigenspaces of the original and the modified network configurations. Third, some lower-dimensional conditions on the discernibility of topological variations are established, which explicitly illustrate how the network topologies, sensor locations, node-system dynamics and output, as well as inner interactions altogether affect the discernibility of the networks. Compared with \cite{haoyuqing} and \cite{indiscernible}, which assume that all state variables of the networks are accessible, these new conditions verify the possibility of detecting topological changes using the output of a subset of nodes, thus are more general and have broader applicability in practice. Fourth, the output discernibility of topological variations for multi-agent systems is revisited, revealing that the sufficiency of the criterion given in \cite{detecting_topology} does not hold, and consequently a modified necessary and sufficient condition is derived.}

The remainder of this paper is organized as follows: Some preliminaries and the model description are given in Section \ref{section2}. Some eigenspace-based conditions on the discernibility of topological changes are developed in Section \ref{section3}. Some lower-dimensional conditions on the discernibility of topological variations are established in Section \ref{section4}. The output discernibility of topological changes for multi-agent systems is revisited in Section \ref{section5}, with a modified condition derived. Finally, conclusions are drawn in Section \ref{section6}.

\section{Preliminaries and model description}
\label{section2}
Some preliminaries and the model description are introduced in this section.
\subsection{Notation and definitions}
Let $\mathbb{N}$, $\mathbb{R}$ and $\mathbb{C}$ be the fields of integers, real and complex numbers, respectively. Let ${I_n}$ be the identity matrix of size $n \times n$, $e_i$ be the vector with all zero entries except for $[e_i]_i=1$, $\textbf{0}_{n }$ be the $n \times 1$ vector with all zero entries, and $diag\left\{ {{a_1},{a_2}, \cdots ,{a_n}} \right\}$ be a diagonal matrix with diagonal entries ${a_1},{a_2}, \cdots ,{a_n}$. The linear span of a set of vectors over the complex field is denoted as $span\{.\}$. Moreover, for a matrix $A \in \mathbb{R}^{n \times n}$, $\sigma (A)=\{\lambda_1,\cdots,\lambda_r\}$ denotes the set of all its eigenvalues, $S(\lambda_i|A)=span\{x \in \mathbb{C}^n|Ax=\lambda_ix\}$ denotes the eigenspace corresponding to $\lambda_i$, and $\tau(\lambda_i |A)$ denotes the geometric multiplicity of $\lambda_i$. The null space of a real matrix $M \in \mathbb{R}^{n \times m}$ is denoted as $\mathcal{N}(M)$. The dimension of a vector space is denoted by \textbf{dim}. Let $A \otimes B$ be the Kronecker product of matrices $A$ and $B$. Given a set of matrices $\{A_1,\cdots,A_n\}$, if they have the same column dimension, then $col(A_1,\cdots,A_n)=[A_1^T,\cdots,A_n^T]^T$. Let ${V_1} + {V_2}$ and ${V_1} \oplus {V_2}$ be the sum and the direct sum of spaces $V_1$ and $V_2$, respectively. Matrices, if their dimensions are not explicitly indicated, are assumed to be compatible for algebraic operations.

{A directed and weighted graph ${\mathcal G}=({\mathcal V},{\mathcal E},{\mathcal W})$ consists of a node set ${\mathcal V}=\{1,\cdots,n\}$, an edge set ${\mathcal E} \subset {\mathcal V} \times {\mathcal V}$, {and a weight matrix ${\mathcal W} =[w_{ij}] \in \mathbb{R}^{n \times n}$}. Note that $(j,i) \in {\mathcal E}$ if and only if $w_{ij} \ne 0$. The adjacency matrix of graph ${\mathcal G}$ is denoted by ${\mathcal A}({\mathcal G})=[a_{ij}]\in \mathbb{R}^{n\times n}$, where $a_{ij}=w_{ij}$ if $(j,i) \in {\mathcal E}${,} and $a_{ij}=0$ otherwise.}

%

\begin{definition}\cite{roman_advanced_2005}
A vector $x_m$ is called an $m$th-order generalized eigenvector of matrix $A$ corresponding to the eigenvalue $\lambda$ if
${(A - \lambda I)^m}{x_m} = 0$ and
${(A - \lambda I)^{m - 1}}{x_m} \ne 0$.
Also, $x_1$, $\cdots$, $x_g$ form a Jordan chain of $A$ with top vector $x_1$, where the maximum number $g$ is called the length of this Jordan chain.

\end{definition}

\begin{definition}
\label{definition4}
\cite{haoyuqing}
Let $A \in \mathbb{C}^{n \times n}$, $H \in \mathbb{C}^{n \times n}$, and $\lambda$ be an eigenvalue of $A$. If vectors $x_1$, $x_2$, $\cdots$, $x_{\theta}$ satisfy $({\lambda }I - A){x _1} = 0$ and $({\lambda}I - A){x _{i + 1}} = H{x _i}$ for $i \in \{ 1, \cdots ,\theta  - 1\}$,
then $x_1$, $x_2$, $\cdots$, $x_{\theta}$ constitute a generalized Jordan chain of $A$ about $H$ corresponding to the eigenvalue $\lambda$, where $x_1$ is the top vector{,} and the maximum number $\theta$ is the length of this generalized Jordan chain.

\end{definition}

\subsection{Model description}
{Consider a network consisting of $N$ identical nodes, with a general directed and weighted topology ${\mathcal G}=({\mathcal V},{\mathcal E},{\mathcal W})$. Each node is represented by an LTI system:
    \begin{equation}
\label{eq1}
    \left\{ \begin{array}{l}
{{\dot x}_i} = A{x_i} + \sum\limits_{j = 1}^N {{w_{ij}}H{x_j}}, \\
{y_i} = C{x_i},
\end{array} \right.\;\;\;\;\;\;i = 1,2, \cdots ,N,
    \end{equation}
   where ${x_i} \in {\mathbb{R}^n}$ is the state vector, {$A \in {\mathbb{R}^{n \times n}}$ is the state matrix describing the dynamics of the node-systems; ${w _{ij}} \in \mathbb{R}$ represents the coupling strength between nodes $i$ and $j$, $H \in {\mathbb{R}^{n \times n}}$ denotes the inner coupling matrix describing the interconnections among components of $x_j$; ${y_i} \in {\mathbb{R}^p}$ is the output vector, and $C \in {\mathbb{R}^{p \times n}}$ is the output matrix. To avoid trivial situations, always assume $N \ge 2$ in this paper. Assume that ${w _{ii}} = 0$, and ${w _{ij}} \ne 0$ if there is an edge from node $j$ to node $i$, otherwise ${w _{ij}} = 0$, for all $i,j = 1,2, \cdots ,N$.}
Let $L = [{w _{ij}}] \in {\mathbb{R}^{N \times N}}$,
which represents the network topology for the networked systems (\ref{eq1}). Let $\Delta=diag\{\delta_1,\;\delta_2,\cdots,\delta_N\}$ be an index matrix indicating the sensor locations, where $\delta_i=1$ if node $i$ is available for measurement, but otherwise $\delta_i=0$, for $i=1,\;2,\cdots,N$.
Moreover, let $X = {[ {x_{\rm{1}}^T,\;x_2^T,\; \cdots ,\;x_N^T} ]^T}$ and $Y = [ y_1^T,\;y_2^T,\; \cdots,y_N^T ]^T$ be the whole state and all output of the networked systems, respectively. Then, the networked systems (\ref{eq1}) can be rewritten in a compact form as
\begin{equation}
\label{eq3}
\left\{
 \begin{array}{l}
 \dot X = \Phi X,\\
 Y=\Psi X,
 \end{array}
 \right.
\end{equation}
where
\begin{equation}
\label{eq4}
 \begin{array}{l}
\Phi  = {I_N} \otimes A + L \otimes H,\;
\Psi = \Delta \otimes C.
 \end{array}
\end{equation}

\section{Eigenspace-based conditions on the $\Psi$-discernibility of topological changes}
\label{section3}
The effect of topological changes on the network output is investigated. In particular, the interest is in characterizing the topological changes that do not alter the output trajectories
(for certain initial states). A topological change caused by a removal/addition of an edge, or a change in an
edge weight, results in a new network
\begin{equation}
\label{eq5}
\left\{\begin{array}{l}
 \dot {\bar X} = {\bar \Phi} {\bar X},\\
 \bar Y=\Psi \bar X,
 \end{array}\right.
\end{equation}
with
\begin{equation}
\label{eq6}
\bar \Phi  = {I_N} \otimes A + \bar L \otimes H.
\end{equation}
Now, the concept of $\Psi$-indiscernible pair of initial states is introduced.

\begin{definition}\cite{detecting_topology}
Consider two networked systems (\ref{eq3})-(\ref{eq4}) and (\ref{eq5})-(\ref{eq6}). A pair of initial states $(X_0,\;\bar X_0) \in \mathbb{R}^{Nn } \times \mathbb{R}^{Nn }$ is called $\Psi$-indiscernible with respect
to the topological change $L \to \bar L$ if and only if
$\Psi e^{\Phi t} X_0=\Psi e^{{\bar \Phi} t} {\bar X}_0,$ for all $t \ge 0$.
\end{definition}

Note that $(X_0,\;\bar X_0)=(\textbf{0}_{Nn},\;\textbf{0}_{Nn})$ is always a $\Psi$-indiscernible pair (irrespective of the specific topological variation), which is called the trivial indiscernible pair.
According to whether a nontrivial $\Psi$-indiscernible pair exists, topological changes are classified into two groups as follows.

\begin{definition}\cite{detecting_topology}
\label{definition2}
For the networked system (\ref{eq3})-(\ref{eq4}), a topological change $L \to \bar L$
is called $\Psi$-discernible if there is no (nontrivial) $\Psi$-indiscernible pair of initial states. Otherwise, it is called $\Psi$-indiscernible.
\end{definition}
%
%
%
%
%
%

%

Let $\mu$ be an eigenvalue of $\Phi$ with the corresponding eigenspace $S(\mu|\Phi)$, and let $\bar \mu$ be an eigenvalue of $\bar \Phi$ with the associated eigenspace $S(\bar \mu|\bar \Phi)$.
In what follows, a necessary and sufficient condition on the $\Psi$-discernibility of topological changes is established in terms of the eigenspaces of the original and the modified networks.

\begin{theorem}
\label{theorem10}
Consider the networked system (\ref{eq3})-(\ref{eq4}). A topological change $L \to \bar L$ is $\Psi$-discernible if and only if for all $\mu \in \sigma(\Phi)\cup \sigma(\bar \Phi)$, the following two conditions hold simultaneously:
\begin{enumerate}
  \item $S(\mu|\Phi) \cap  S(\mu|\bar \Phi)=\{\textbf{0}_{Nn}\}$;
  \item $\mathcal {N}(\Psi) \cap \{S(\mu|\Phi) \oplus  S(\mu|\bar\Phi)\}=\{\textbf{0}_{Nn}\}$.
\end{enumerate}

\end{theorem}

{\textbf{Proof:}}
It is easy to verify that the topological change $L \to \bar L$ is $\Psi$-discernible if and only if $(diag \{\Phi,\;\bar\Phi\}, [
\Psi \; -\Psi])$ is observable.

Necessity: If there exists $\mu^* \in \sigma(\Phi)\cup \sigma(\bar \Phi)$ such that $S(\mu^*|\Phi) \cap S(\mu^*|\bar \Phi)\ne\{\textbf{0}_{Nn}\}$, then there exists a nonzero vector $x \in S(\mu^*|\Phi) \cap S(\mu^*|\bar \Phi)$. It is easy to verify that $\Phi x=\mu^* x$ and $\bar \Phi x=\mu^* x$. Thus, there exists an eigenpair $(\mu^*,[ {\begin{array}{*{20}{c}}
x^T& x^T
\end{array}} ]^T)$ satisfying $\left[ {\begin{array}{*{20}{c}}
\Phi&\\
&\bar \Phi
\end{array}} \right]\left[ {\begin{array}{*{20}{c}}
x\\
x
\end{array}} \right]=\mu^*\left[ {\begin{array}{*{20}{c}}
x\\
x
\end{array}} \right]$ and
$[ {\begin{array}{*{20}{c}}
\Psi&-\Psi
\end{array}} ]\left[ {\begin{array}{*{20}{c}}
x\\
x
\end{array}} \right]=\textbf{0}$, which implies that $(diag \{\Phi,\bar\Phi\}, [\Psi\; -\Psi])$ is unobservable. Therefore, the topological change $L \to \bar L$ is $\Psi$-indiscernible.

Assume that for all $\mu \in \sigma(\Phi) \cup \sigma(\bar \Phi)$, $S(\mu|\Phi) \cap S(\mu|\bar \Phi)=\{\textbf{0}_{Nn}\}$. If there exists $\mu^* \in \sigma(\Phi)\cup \sigma(\bar \Phi)$ such that $\mathcal{N} (\Psi) \cap \{S(\mu^*|\Phi) \oplus S(\mu^*|\bar\Phi)\} \ne \{\textbf{0}_{Nn}\}$, then there exists a nonzero vector $x \in \{S(\mu^*|\Phi) \oplus S(\mu^*|\bar\Phi)\}$ such that $\Psi x =\textbf{0}$. Let $x=x_1+x_2$, where $x_1 \in S(\mu^*|\Phi)$ and $x_2 \in S(\mu^*|\bar\Phi)$. From $x \ne \textbf{0}_{Nn}$, it follows that $x_1$ and $x_2$ are not both $\textbf{0}_{Nn}$, as discussed in the following three cases:

\begin{itemize}
  \item If $x_2=\textbf{0}_{Nn}$ and $x_1 \ne \textbf{0}_{Nn}$, it is easy to verify that $\Phi x_1=\mu^*x_1 $ and $\Psi x_1=0$. Then, there exists a nonzero vector $[ {\begin{array}{*{20}{c}}
x_1^T& \textbf{0}_{Nn}^T
\end{array}} ]^T $ satisfying $\left[ {\begin{array}{*{20}{c}}
\Phi&\\
&\bar \Phi
\end{array}} \right]\left[ {\begin{array}{*{20}{c}}
x_1\\
\textbf{0}_{Nn}
\end{array}} \right]=\mu^*\left[ {\begin{array}{*{20}{c}}
x_1\\
\textbf{0}_{Nn}
\end{array}} \right]$ and
$[ {\begin{array}{*{20}{c}}
\Psi&-\Psi
\end{array}} ]\left[ {\begin{array}{*{20}{c}}
x_1\\
\textbf{0}_{Nn}
\end{array}} \right]=\textbf{0}$, which implies that $(diag \{\Phi,\bar\Phi\}, [\Psi\; -\Psi])$ is unobservable. Thus, the topological change $L \to \bar L$ is $\Psi$-indiscernible.
  \item If $x_1=\textbf{0}_{Nn}$ and $x_2 \ne \textbf{0}_{Nn}$, one can similarly prove that the topological change $L \to \bar L$ is $\Psi$-indiscernible.
  \item If $x_1 \ne \textbf{0}_{Nn}$ and $x_2 \ne \textbf{0}_{Nn}$, then $\Phi x_1=\mu^* x_1$, $\bar \Phi x_2=\mu^* x_2$ and $\Psi (x_1+x_2)=\textbf{0}$. Consequently, there exists a nonzero vector $[ {\begin{array}{*{20}{c}}
x_1^T& -x_2^T
\end{array}} ]^T $ satisfying $\left[ {\begin{array}{*{20}{c}}
\Phi&\\
&\bar \Phi
\end{array}} \right]\left[ {\begin{array}{*{20}{c}}
x_1\\
-x_2
\end{array}} \right]=\mu^*\left[ {\begin{array}{*{20}{c}}
x_1\\
-x_2
\end{array}} \right]$ and
$[ {\begin{array}{*{20}{c}}
\Psi&-\Psi
\end{array}} ]\left[ {\begin{array}{*{20}{c}}
x_1\\
-x_2
\end{array}} \right]=\textbf{0}$, which implies that $(diag \{\Phi,\bar\Phi\}, [\Psi\; -\Psi])$ is unobservable. Thus, the topological change $L \to \bar L$ is $\Psi$-indiscernible.
\end{itemize}
Therefore, if there exists $\mu^* \in \sigma(\Phi)\cup \sigma(\bar \Phi)$ such that $\mathcal{N} (\Psi) \cap \{S(\mu^*|\Phi) \oplus S(\mu^*|\bar\Phi)\} \ne \{\textbf{0}_{Nn}\}$, then the topological change $L \to \bar L$ is $\Psi$-indiscernible.

Sufficiency: If the topological change $L \to \bar L$ is $\Psi$-indiscernible, then there exists an eigenpair of $diag\{\Phi,\;\bar \Phi\}$, denoted as $(\mu^*,\eta)$, such that $[ {\begin{array}{*{20}{c}}
\Psi&-\Psi
\end{array}} ] \eta=\textbf{0}$. Let $\textbf{0}\ne \eta=[ {\begin{array}{*{20}{c}}
\eta_1^T& \eta_2^T
\end{array}} ]^T$, where $\eta_1,\;\eta_2 \in \mathbb{C}^{Nn}$. Then, $\eta_1$ and $\eta_2$ are not both $\textbf{0}_{Nn}$, and they satisfy \[\left\{ \begin{array}{l}
\Phi \eta_1 = \mu^* \eta_1,\\
\bar \Phi \eta_2 =\mu^* \eta_2,\\
\Psi(\eta_1-\eta_2)=\textbf{0}.
\end{array} \right.\]
It follows that $\eta_1 \in S(\mu^*|\Phi)$ and $-\eta_2 \in  S (\mu^*|\bar\Phi)$, thus $\eta_1-\eta_2 \in \{S(\mu^*|\Phi) +  S(\mu^*|\bar\Phi)\}$. If $S(\mu^*|\Phi) \cap S(\mu^*|\bar\Phi)\ne\{\textbf{0}_{Nn}\}$, then condition $(1)$ does not hold. If $S(\mu^*|\Phi) \cap S(\mu^*|\bar\Phi)=\{\textbf{0}_{Nn}\}$, then $\eta_1-\eta_2 \in \{S(\mu^*|\Phi) \oplus  S(\mu^*|\bar\Phi)\}$. Noting that $\eta_1$ and $\eta_2$ are not both $\textbf{0}_{Nn}$, one can conclude that $\eta_1-\eta_2 \ne \textbf{0}_{Nn}$. Thus, $\textbf{0}_{Nn} \ne \eta_1-\eta_2 \in \mathcal{N} (\Psi) \cap \{S(\mu^*|\Phi) \oplus  S(\mu^*|\bar\Phi)\}$, which implies that condition $(2)$ does not hold. Therefore, if the topological change $L \to \bar L$ is $\Psi$-indiscernible, then at least one condition in Theorem \ref{theorem10} does not hold.
\hfill $\blacksquare$

\begin{remark}
The possibility of detecting topological variations by observing network states was investigated in \cite{indiscernible}, where it was claimed that the topological variation is always-discernible if and only if $\Phi$ and $\bar \Phi$ have no common eigenpairs. This condition can be reproduced from Theorem \ref{theorem10} as a special case with $\Psi=I$. Theorem \ref{theorem10} extends the results in \cite{indiscernible} to a more general and practical case of detecting topological changes from the observation of output trajectories. This nontrivial extension is of great significance for engineering applications.

\end{remark}

\begin{remark}
Noting that $S(\mu|\Phi) \subseteq  \{S(\mu|\Phi)\oplus  S(\mu|\bar\Phi)\}$, the second condition in Theorem \ref{theorem10} requires that $\mathcal {N}(\Psi) \cap S(\mu|\Phi) =\{\textbf{0}_{Nn}\}$ for all $\mu \in \sigma(\Phi)$, which implies that $(\Phi, \Psi)$ is observable. Similarly, since $S(\mu|\bar \Phi) \subseteq  \{S(\mu|\Phi)\oplus  S(\mu|\bar\Phi)\}$, Theorem \ref{theorem10} requires the observability of $(\bar\Phi, \Psi)$ as well. Therefore, the observability of $(\Phi, \Psi)$ and $(\bar\Phi, \Psi)$ is necessary for the $\Psi$-discernibility of the topological variation. Moreover, if $\Phi$ and $\bar \Phi$ have no common eigenvalues, then the topological change $L \to \bar L$ is $\Psi$-discernible if and only if both $(\Psi, \Phi)$ and $(\Psi, \bar \Phi)$ are observable.

\end{remark}

Let $f:\;\mathbb{C}^{Nn} \to \mathbb{C}^{Np}$ be a linear map such that, for $x \in \mathbb{C}^{Nn}$, $f(x)=\Psi x$. Moreover, let $\tau(\mu|\Phi)$ and $\tau(\mu|\bar\Phi)$ denote the geometric multiplicities of $\mu$ for $\Phi$ and $\bar \Phi$, respectively. In what follows, some conditions on the $\Psi$-discernibility of topological variations are derived in terms of the multiplicities of the eigenvalues.

\begin{corollary}
\label{theorem20}
Consider the networked system (\ref{eq3})-(\ref{eq4}). A topological change $L \to \bar L$ is $\Psi$-discernible if and only if $\textbf{dim} \{f[S(\mu|\Phi)]+f[ S(\mu|\bar\Phi)]\}=\tau(\mu|\Phi)+ \tau(\mu|\bar\Phi)$, for all $\mu \in \sigma(\Phi)\cup \sigma(\bar \Phi)$.

\end{corollary}

\textbf{Proof:}
Note that $\textbf{dim} \{f[S(\mu|\Phi)]+f[ S(\mu|\bar\Phi)]\}=\textbf{dim}\{f[S(\mu|\Phi)]\}+\textbf{dim}\{f[ S(\mu|\bar\Phi)]\}-\textbf{dim}\{f[S(\mu|\Phi)]\cap f[ S(\mu|\bar\Phi)]\}=\textbf{dim}\{S(\mu|\Phi)\}-\textbf{dim}\{\mathcal N(\Psi) \cap S(\mu|\Phi)\}+\textbf{dim}\{ S(\mu|\bar\Phi)\}-\textbf{dim}\{\mathcal N(\Psi) \cap S(\mu|\bar \Phi)\}-\textbf{dim}\{f[S(\mu|\Phi)]\cap f[ S(\mu|\bar\Phi)]\}=\tau(\mu|\Phi)+ \tau(\mu|\bar\Phi)-\textbf{dim}\{\mathcal N(\Psi) \cap S(\mu|\Phi)\}-\textbf{dim}\{\mathcal N(\Psi) \cap S(\mu|\bar \Phi)\}-\textbf{dim}\{f[S(\mu|\Phi)]\cap f[ S(\mu|\bar\Phi)]\}$. Thus, $\textbf{dim} \{f[S(\mu|\Phi)]+f[ S(\mu|\bar\Phi)]\}=\tau(\mu|\Phi)+ \tau(\mu|\bar\Phi)$ if and only if the following three conditions hold simultaneously:
\begin{equation}
\label{eq100}
\mathcal N(\Psi) \cap S(\mu|\Phi)=\{\textbf{0}_{Nn}\};
\end{equation}
\begin{equation}
\label{eq101}
\mathcal N(\Psi) \cap S(\mu|\bar \Phi)=\{\textbf{0}_{Nn}\};
\end{equation}
\begin{equation}
\label{eq102}
f[S(\mu|\Phi)]\cap f[ S(\mu|\bar\Phi)]=\{\textbf{0}_{Nn}\}.
\end{equation}

Necessity: If there exists $\mu^* \in \sigma(\Phi)\cup \sigma(\bar \Phi)$ such that $\textbf{dim} \{f[S(\mu^*|\Phi)]+f[ S(\mu^*|\bar\Phi)]\}<\tau(\mu^*|\Phi)+ \tau(\mu^*|\bar\Phi)$, then at least one of (\ref{eq100}), (\ref{eq101}), (\ref{eq102}) does not hold.
\begin{itemize}
  \item Case 1: If (\ref{eq100}) does not hold, then there exists a nonzero vector $v \in S(\mu^*|\Phi)$ such that $\Psi v=\textbf{0}$, which implies that $\textbf{0} \ne v \in \mathcal {N}(\Psi) \cap \{S(\mu^*|\Phi) \oplus  S(\mu^*|\bar\Phi)\}$. Thus, it follows from Theorem \ref{theorem10} that the topological change $L \to \bar L$ is $\Psi$-indiscernible.
  \item Case 2: If (\ref{eq101}) does not hold, one can similarly prove that the topological change $L \to \bar L$ is $\Psi$-indiscernible.
  \item Case 3: If (\ref{eq102}) does not hold,
  then there exist $v \in S(\mu^*|\Phi)$ and $w \in S(\mu^*|\bar \Phi)$ such that $\Psi v=\Psi w \ne \textbf{0}$. If $S(\mu^*|\Phi) \cap S(\mu^*|\bar\Phi) \ne \{\textbf{0}_{Nn}\}$, it follows from Theorem \ref{theorem10} that the topological change $L \to \bar L$ is $\Psi$-indiscernible; If $S(\mu^*|\Phi) \cap S(\mu^*|\bar\Phi) = \{\textbf{0}_{Nn}\}$, then $v-w \in \{S(\mu^*|\Phi) \oplus S(\mu^*|\bar\Phi)\}$. From $\Psi v=\Psi w \ne \textbf{0}$, it follows that $v \ne \textbf{0}$ and $w \ne \textbf{0}$. Then, one can show that $v-w \ne \textbf{0} $. Since $\Psi (v-w)=\textbf{0}$, one gets that $\textbf{0} \ne v-w \in \mathcal N (\Psi) \cap \{S(\mu^*|\Phi) \oplus S(\mu^*|\bar\Phi)\}$. According to Theorem \ref{theorem10}, the topological change $L \to \bar L$ is $\Psi$-indiscernible.
  \end{itemize}
Therefore, if there exists $\mu^* \in \sigma(\Phi)\cup \sigma(\bar \Phi)$ such that $\textbf{dim} \{f[S(\mu^*|\Phi)]+f[ S(\mu^*|\bar\Phi)]\}<\tau(\mu^*|\Phi)+ \tau(\mu^*|\bar\Phi)$, then the topological change $L \to \bar L$ is $\Psi$-indiscernible.

Sufficiency: If the topological change $L \to \bar L$ is $\Psi$-indiscernible, then at least one condition in Theorem \ref{theorem10} does not hold.
\begin{itemize}
  \item Case 1: There exists $\mu^* \in \sigma(\Phi) \cup \sigma(\bar\Phi)$ such that $S(\mu^*|\Phi) \cap S(\mu^*|\bar\Phi) \ne \{\textbf{0}_{Nn}\}$. Then, there exists a nonzero vector $x \in S(\mu^*|\Phi) \cap S(\mu^*|\bar\Phi)$. This will be discussed in the following two cases.
      \begin{description}
        \item[-] If $x \notin \mathcal N(\Psi)$, then $\textbf{0} \ne f(x) \in f[S(\mu^*|\Phi)]\cap f[ S(\mu^*|\bar\Phi)$, which indicates that (\ref{eq102}) does not hold. Thus, $\textbf{dim} \{f[S(\mu^*|\Phi)]+f[ S(\mu^*|\bar\Phi)]\}<\tau(\mu^*|\Phi)+ \tau(\mu^*|\bar\Phi)$.
        \item[-] If $x \in \mathcal N(\Psi)$, then $\textbf{0} \ne x \in \mathcal N(\Psi) \cap S(\mu^*|\Phi)$, which implies that (\ref{eq100}) does not hold. Thus, $\textbf{dim} \{f[S(\mu^*|\Phi)]+f[ S(\mu^*|\bar\Phi)]\}<\tau(\mu^*|\Phi)+ \tau(\mu^*|\bar\Phi)$.
      \end{description}

  \item Case 2: There exists $\mu^* \in \sigma(\Phi) \cup \sigma(\bar\Phi)$ such that $\mathcal N(\Psi) \cap \{S(\mu^*|\Phi) \oplus S(\mu^*|\bar \Phi)\}\ne \{\textbf{0}_{Nn}\}$. Then, there exists a nonzero vector $x \in S(\mu^*|\Phi) \oplus S(\mu^*|\bar \Phi)$ such that $\Psi x=f(x)=\textbf{0}$. Let $x=x_1+x_2$, where $x_1 \in S(\mu^*|\Phi)$ and $x_2 \in S(\mu^*|\bar \Phi)$. Then, it follows from $f(x)=f(x_1+x_2)=f(x_1)+f(x_2)=\textbf{0}$ that $f(x_1)=f(-x_2)$. Note that $x_1 \in S(\mu^*|\Phi)$ and $-x_2 \in S(\mu^*|\bar \Phi)$.
       If $f(x_1)=f(-x_2) \ne \textbf{0}$, then (\ref{eq102}) does not hold. Thus, $\textbf{dim} \{f[S(\mu^*|\Phi)]+f[ S(\mu^*|\bar\Phi)]\}<\tau(\mu^*|\Phi)+ \tau(\mu^*|\bar\Phi)$. In what follows, the case of $f(x_1)=f(-x_2)=\textbf{0}$ will be discussed. Noting that $x=x_1+x_2 \ne \textbf{0}$, one can verify that $x_1$ and $x_2$ are not both $\textbf{0}$.
      If $x_1 \ne \textbf{0}$, then $\textbf{0} \ne x_1 \in \mathcal N (\Psi) \cup S(\mu^*|\Phi)$. Thus, (\ref{eq100}) does not hold, which implies that $\textbf{dim} \{f[S(\mu^*|\Phi)]+f[ S(\mu^*|\bar\Phi)]\}<\tau(\mu^*|\Phi)+ \tau(\mu^*|\bar\Phi)$.
      If $x_1=\textbf{0}$, then $x_2 \ne \textbf{0}$. It follows that $\textbf{0} \ne -x_2 \in \mathcal N (\Psi) \cup S(\mu^*|\bar \Phi)$, which indicates that (\ref{eq101}) does not hold. Thus, $\textbf{dim} \{f[S(\mu^*|\Phi)]+f[ S(\mu^*|\bar\Phi)]\}<\tau(\mu^*|\Phi)+ \tau(\mu^*|\bar\Phi)$.
      \end{itemize}

 Therefore, if the topological change $L \to \bar L$ is $\Psi$-indiscernible, then there exists $\mu^* \in \sigma(\Phi) \cup \sigma(\bar\Phi)$ such that $\textbf{dim }\{f[S(\mu^*|\Phi)]+f[ S(\mu^*|\bar\Phi)]\} < \tau(\mu^*|\Phi)+ \tau(\mu^*|\bar\Phi)$.
\hfill $\blacksquare$

In what follows, a necessary condition on the $\Psi$-discernibility is derived, which is intuitive and easier to verify.

\begin{corollary}
\label{corollary4}
Consider the networked system (\ref{eq3})-(\ref{eq4}). If the topological change $L \to \bar L$ is $\Psi$-discernible, then $Rank(\Delta) \times Rank(C)\ge \tau(\mu|\Phi)+ \tau(\mu|\bar\Phi)$, for all $\mu \in \sigma(\Phi)\cup \sigma(\bar \Phi)$.

\end{corollary}

\textbf{Proof:}
Assume that there exists $\mu^* \in \sigma(\Phi)\cup \sigma(\bar \Phi)$ such that $Rank(\Delta) \times Rank(C)<\tau(\mu^*|\Phi)+ \tau(\mu^*|\bar\Phi)$. If $S(\mu^*|\Phi) \cap  S(\mu^*|\bar\Phi) \ne \{\textbf{0}_{Nn}\}$, it follows from Theorem \ref{theorem10} that the topological change $L \to \bar L$ is $\Psi$-indiscernible. If $S(\mu^*|\Phi) \cap  S(\mu^*|\bar\Phi) = \{\textbf{0}_{Nn}\}$, then $\textbf{dim}\{S(\mu^*|\Phi) \oplus  S(\mu^*|\bar\Phi)\}=\textbf{dim}[S(\mu^*|\Phi)]+\textbf{dim}[ S(\mu^*|\bar\Phi)]=\tau(\mu^*|\Phi)+ \tau(\mu^*|\bar\Phi)>Rank(\Delta) \times Rank(C)$. Noting that $\textbf{dim}\{\mathcal N(\Psi)\}+Rank(\Delta) \times Rank(C)=nN$, it can be easily verified that $\textbf{dim}\{\mathcal N(\Psi)\}\\+\textbf{dim}\{S(\mu^*|\Phi) \oplus  S(\mu^*|\bar\Phi)\}>Nn$. It follows that $\mathcal N(\Psi) \cap \{S(\mu^*|\Phi) \oplus  S(\mu^*|\bar\Phi)\} \ne \{\textbf{0}_{Nn}\}$. According to Theorem \ref{theorem10}, the topological change $L \to \bar L$ is $\Psi$-indiscernible. Therefore, if there exists $\mu^* \in \sigma(\Phi)\cup \sigma(\bar \Phi)$ such that $Rank(\Delta) \times Rank(C)<\tau(\mu^*|\Phi)+ \tau(\mu^*|\bar\Phi)$, then the topological change $L \to \bar L$ is $\Psi$-indiscernible.

\hfill $\blacksquare$

The effectiveness of the above corollary is illustrated via the following example.

\begin{example}
Consider a simple network consisting of four connected identical nodes, shown in (a) of Fig. \ref{fig10}, with $w_{21}=w_{32}=w_{43}=w_{14}=1$. Suppose that the output of the first node and the second node can be observed, i.e., $\delta_1=\delta_2=1$, with
\[A = \left[ {\begin{array}{*{20}{c}}
1&1\\
0&2
\end{array}} \right],\;B = \left[ {\begin{array}{*{20}{c}}
1&0\\
0&0
\end{array}} \right],\;C = \left[ {\begin{array}{*{20}{c}}
1&0
\end{array}} \right].\]
\begin{figure}[htb]
\centering
\subfloat[${\mathcal G}$]{%
\includegraphics[width=.15\textwidth]{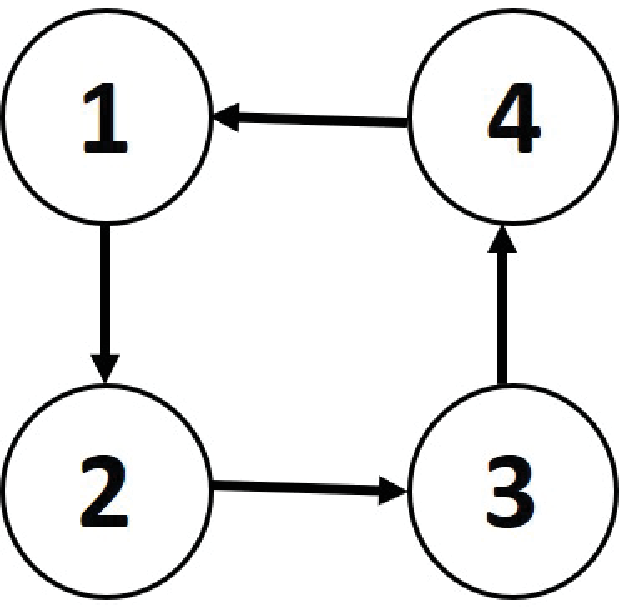}}
\hspace{.25in}
\subfloat[$\bar {\mathcal G}$]{%
\includegraphics[width=.15\textwidth]{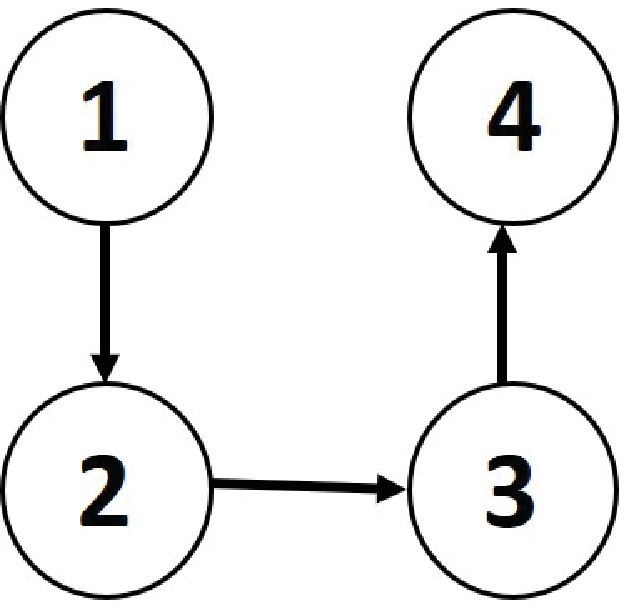}}\hfill
\caption{Network topologies}
\label{fig10}
\end{figure}

It can be easily verified that $L = \left[ {\begin{array}{*{20}{c}}
0&{ 0}&{ 0}&1\\
{ 1}&0&{ 0}&0\\
{ 0}&{ 1}&0 &0\\
0&0&1&0
\end{array}} \right]$ and $\Delta=\left[ {\begin{array}{*{20}{c}}
{ 1}&0&{ 0}&0\\
{ 0}&{ 1}&0&0
\end{array}} \right]$. Assume that the edge from node $4$ to node $1$ is removed. Then, the new network is shown in (b) of Fig. \ref{fig10}, with the topology matrix $\bar L = \left[ {\begin{array}{*{20}{c}}
0&{ 0}&0&0\\
{ 1}&0&{0}&0\\
0&{ 1}&0&0\\
0&0&1&0
\end{array}} \right]$.

It is easy to verify that $\sigma(\Phi)=\{0,\;2,\;1+i,\;1-i\}$ and $\tau(2|\Phi)=4$. Moreover, $\sigma(\bar \Phi)=\{1,\;2\}$ and $\tau(2|\bar\Phi)=4$. It follows that $\tau(2|\Phi)+\tau(2|\bar\Phi)=8$. Thus, $Rank(\Delta) \times Rank(C)=2<\tau(2|\Phi)+\tau(2|\bar\Phi)$. According to Corollary \ref{corollary4}, the topological change $L \to \bar L$ is $\Psi$-indiscernible.

\end{example}

\begin{remark}
It follows from Corollary \ref{corollary4} that to guarantee the $\Psi$-discernibility of the topological variation, the number of the sensors should not be less than $\frac{\mathop  {max} \limits_{\mu \in \sigma(\Phi) \cup \sigma(\bar \Phi)}{\tau(\mu|\Phi)+\tau(\mu|\bar\Phi)}}{Rank C}$. Thus, Corollary \ref{corollary4} provides useful information on the number of the sensors required for the $\Psi$-discernibility.

\end{remark}

\section{Lower-dimensional conditions on the $\Psi$-discernibility of topological changes}
\label{section4}
In Section \ref{section3}, some eigenspace-based conditions on
the $\Psi$-discernibility of topological changes are developed. In this section, some
lower-dimensional conditions are further established by taking the network structures into consideration, which are easier to verify and apply.

Let $\lambda_1$, $\lambda_2$, $\cdots$, $\lambda_s$ be the eigenvalues of $L$. The Jordan chain of $L$ associated with the eigenvalue $\lambda_i$ is denoted as $t_i^1$, $t_i^2$, $\cdots$, $t_i^{\alpha_i}$, where $t_i^1$ is the top vector, and $\alpha_i$ is the length of the Jordan chain. Moreover, denote the eigenvalues of $A+\lambda_iH$ as $\mu_i^j$, with the corresponding generalized Jordan chain about $H$ denoted as $\xi_{ij}^1$, $\xi_{ij}^2$, $\cdots$, $\xi_{ij}^{\theta_{ij}}$, where $\xi_{ij}^1$ is the top vector, and $\theta_{ij}$ is the length, for $j=1,\cdots,p_i$, $i=1,\cdots,s$. In \cite{haoyuqing}, the eigenspaces of $\Phi$ are expressed through the generalized eigenvectors of some matrices with lower dimensions.

\begin{lemma}
\cite{haoyuqing}\label{lemma1}
Let $\lambda_1$, $\lambda_2$, $\cdots$, $\lambda_s$ be the eigenvalues of $L$, and $\mu_i^1$, ${\cdots}$, $\mu_i^{p_i}$ be the eigenvalues of $A+\lambda_i H$, counting geometric multiplicities. Then, $\mu_1^1$, ${\cdots}$, $\mu_1^{p_1}$, ${\cdots}$, $\mu_s^1$, ${\cdots}$, $\mu_s^{p_s}$ are the eigenvalues of $\Phi$. Moreover, the eigenspace of $\Phi$ corresponding to $\mu_i^j$ is ${V_{ij}} = span\{ \eta _{ij}^1,\eta _{ij}^2, \cdots ,\eta _{ij}^{{\gamma _{ij}}}\} $, where $\eta_{ij}^1=t_i^1 \otimes \xi_{ij}^1$, $\eta_{ij}^2=t_i^2 \otimes \xi_{ij}^1+t_i^1 \otimes \xi_{ij}^2$, $\cdots$, $\eta_{ij}^{\gamma_{ij}}=t_i^{\gamma_{ij}} \otimes \xi_{ij}^1+t_i^{\gamma_{ij}-1} \otimes \xi_{ij}^2+\cdots+t_i^1 \otimes \xi_{ij}^{\gamma_{ij}}$, $\gamma_{ij}=\min \left\{ {{\alpha _i},{\theta _{ij}}} \right\}$, $j=1,\cdots,p_i$, $i=1,\cdots,s$. Let $\Gamma(\mu) =\{(i,j) \in \mathbb{N} \times \mathbb{N}| \mu_i^j=\mu,\;1 \le j\le p_i,\;1 \le i \le s\}$. Then, $S(\mu|\Phi)={\mathop  \oplus \limits_{(i,j)\in \Gamma(\mu)} {{ V}_{{i}{j}}}}$.
\end{lemma}

Similarly, the eigenspaces of $\bar \Phi$ can be expressed through the generalized eigenvectors of some matrices with lower dimensions. Let $\bar \lambda_1$, $\bar \lambda_2$, $\cdots$, $\bar \lambda_{\bar s}$ be the eigenvalues of $\bar L$, and $\bar \mu_i^1$, ${\cdots}$, $\bar \mu_i^{\bar p_i}$ be the eigenvalues of $A+\bar \lambda_i H$, counting geometric multiplicities. Then, $\bar \mu_1^1$, ${\cdots}$, $\bar \mu_1^{\bar p_1}$, ${\cdots}$, $\bar \mu_{\bar s}^1$, ${\cdots}$, $\bar \mu_{\bar s}^{\bar p_{\bar s}}$ are the eigenvalues of $\bar \Phi$. Moreover, the eigenspace of $\bar \Phi$ corresponding to $\bar \mu_i^j$ is ${\bar V_{ij}} = span\{ \bar \eta _{ij}^1,\bar\eta _{ij}^2, \cdots ,\bar \eta _{ij}^{{\bar \gamma _{ij}}}\} $, where $\bar \eta_{ij}^1$, $\cdots$, $\bar \eta_{ij}^{\bar \gamma_{ij}}$ can be expressed using the generalized eigenvectors of $\bar L$ and $A+\bar \lambda_i H$ as shown in Lemma \ref{lemma1}, $j=1,\cdots,\bar p_i$, $i=1,\cdots,\bar s$. Let $\bar \Gamma(\mu) =\{(i,j) \in \mathbb{N} \times \mathbb{N}| \bar \mu_i^j=\mu,\;1 \le j\le \bar p_i,\;1 \le i \le \bar s\}$. Then, $S(\mu|\bar\Phi)={\mathop  \oplus \limits_{(i,j)\in \bar \Gamma(\mu)} {{ \bar V}_{{i}{j}}}}$.

Next, a necessary and sufficient condition on the $\Psi$-discernibility of topological changes is established.

\begin{theorem}
\label{theorem1}
Consider the networked system (\ref{eq3})-(\ref{eq4}). A topological change $L \to \bar L$ is $\Psi$-discernible if and only if for all $\mu \in \sigma(\Phi)\cup \sigma(\bar \Phi)$, the following two conditions hold simultaneously:
\begin{enumerate}
  \item $\left\{{\mathop  \oplus \limits_{(i,j)\in  \Gamma(\mu)} {{  V}_{{i}{j}}}}\right\} \bigcap \left\{{\mathop  \oplus \limits_{(i,j)\in \bar \Gamma(\mu)} {{ \bar V}_{{i}{j}}}}\right\}=\{\textbf{0}_{Nn}\}$;
  \item $\mathcal N (\Delta \otimes C) \cap \left\{{\left\{{\mathop  \oplus \limits_{(i,j)\in  \Gamma(\mu)} {{  V}_{{i}{j}}}}\right\}  \bigoplus \left\{{{\mathop  \oplus \limits_{(i,j)\in \bar \Gamma(\mu)} {{ \bar V}_{{i}{j}}}}}\right\}}\right\}=\{\textbf{0}_{Nn}\}$.
\end{enumerate}

\end{theorem}

\textbf{Proof:}
The proof follows from Theorem \ref{theorem10} and Lemma \ref{lemma1} directly, thus is omitted.
\hfill $\blacksquare$

The effectiveness of Theorem \ref{theorem1} is demonstrated by the following example.

\begin{example}
Consider a simple network of three connected identical nodes, shown in (a) of Fig. \ref{fig1}, with $w_{21}=w_{32}=w_{23}=1$. Suppose that the output of the first node and the third node can be observed, i.e., $\delta_1=\delta_3=1$, with
\[A = \left[ {\begin{array}{*{20}{c}}
1&0\\
1&1
\end{array}} \right], \;H = \left[ {\begin{array}{*{20}{c}}
1&0\\
0&1
\end{array}} \right],\;C=\left[ {\begin{array}{*{20}{c}}
1&0\\
0&1
\end{array}} \right].\]

It can be easily verified that $L = \left[ {\begin{array}{*{20}{c}}
0&0&0\\
1&0&1\\
0&1&0
\end{array}} \right]$ and $\Delta = \left[ {\begin{array}{*{20}{c}}
1&0&0\\
0&0&1
\end{array}} \right]$. The eigenvalues of $L$ are $\lambda_1=0$, $\lambda_2=1$ and $\lambda_3=-1$, with the corresponding eigenvectors $t_1=e_1-e_3$, $t_2=e_2+e_3$ and $t_3=e_2-e_3$, respectively. Then, the eigenvalue of $A+\lambda_1H=A$ is $\mu_1^1=1$, with the corresponding eigenvector $\xi_{11}^1=e_2$, and the eigenvalue of $A+\lambda_2 H=A+H$ is $\mu_2^1=2$, with the corresponding eigenvector $\xi_{21}^1=e_2$. The eigenvalue of $A+\lambda_3 H=A-H$ is $\mu_3^1=0$, with the corresponding eigenvector $\xi_{31}^1=e_2$.
Thus, it follows that $S(1|\Phi)=V_{11}=span\{(e_1-e_3)\otimes e_2\}$, $S(2|\Phi)=V_{21}=span\{(e_2+e_3)\otimes e_2\}$, and $S(0|\Phi)=V_{31}=span\{(e_2-e_3)\otimes e_2\}$.

Assume that the edge from node $3$ to node $2$ is removed. Then, the new topology matrix is $\bar L = \left[ {\begin{array}{*{20}{c}}
0&0&0\\
1&0&0\\
0&1&0
\end{array}} \right]$. The eigenvalue of $\bar L$ is $\bar \lambda_1=0$, with the corresponding Jordan chain $\bar t_1^1=e_3$, $\bar t_1^2=e_2$, $\bar t_1^3=e_1$. So, the eigenvalue of $A+\bar \lambda_1H=A$ is $\bar \mu_1^1=1$, with the generalized Jordan chain $\bar \xi_{11}^1=e_2$, $\bar \xi_{11}^2=-e_1$.
Thus, it follows that $S(1|\bar\Phi)=\bar V_{11}=span\{e_3\otimes e_2,\;e_2\otimes e_2+e_3 \otimes (-e_1)\}$.

Noting that $ S(0|\bar\Phi) = S(2|\bar\Phi)=\{\textbf{0}\}$, one can easily verify that $S(0|\Phi) \cap  S(0|\bar\Phi) =S(2|\Phi) \cap S(2|\bar\Phi) =S(1|\Phi) \cap  S(1|\bar\Phi) =\{\textbf{0}\}$. Moreover, $\mathcal N(\Delta \otimes C) \cap S(0|\Phi)=\mathcal N (\Delta \otimes C) \cap S(2|\Phi)=\mathcal N (\Delta \otimes C) \cap \{S(1|\Phi) \oplus  S(1|\bar\Phi)\}=\{\textbf{0}\}$. Therefore, it follows from Theorem \ref{theorem1} that the topological change $L \to \bar L$ is $\Psi$-discernible.

\begin{figure}[h]
\centering
\subfloat[${\mathcal G}$]{%
\begin{tabular}{c}
\includegraphics[width=.2\textwidth]{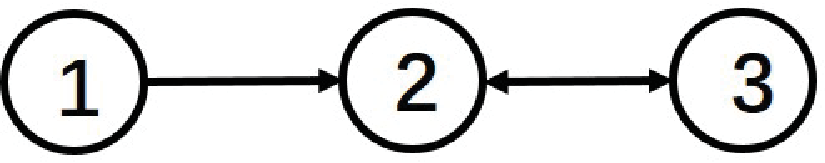}
\end{tabular}
}%
\subfloat[$\bar {\mathcal G}$]{%
\begin{tabular}{c}
\includegraphics[width=.2\textwidth]{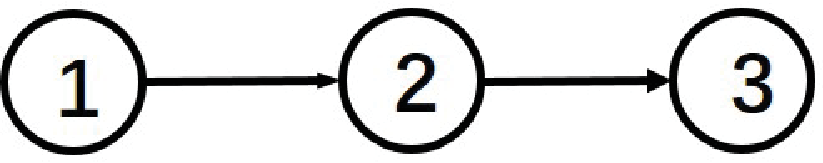}
\end{tabular}
}%
\caption{Network topologies}
\label{fig1}
\end{figure}

\end{example}

In the following, some lower-dimensional and easily-verified conditions on the $\Psi$-discernibility of topological variations are presented.

\begin{corollary}
\label{corollary6}
If a topological change $L \to \bar L$ for the networked system (\ref{eq3})-(\ref{eq4}) is $\Psi$-discernible, then the topological change $L \to \bar L$ is $\Delta$-discernible for the system
\begin{equation}
\label{eq11}
\left\{
 \begin{array}{l}
 \dot x = L x,\\
 y=\Delta x.
 \end{array}
 \right.
 \end{equation}
\end{corollary}

\textbf{Proof:}
Let $\lambda$ be an eigenvalue of $L$, with the corresponding eigenspace
$S(\lambda|L)$.
Moreover, let $\bar\lambda$ be an eigenvalue of $\bar L$, with the associated eigenspace $S(\bar\lambda|\bar L)$.
If the topological change $L \to \bar L$ is $\Delta$-indiscernible for system (\ref{eq11}), then at least one of the following cases occurs.
\begin{itemize}
  \item There exists $\lambda^* \in \sigma(L) \cup \sigma(\bar L)$ such that $S(\lambda^*|L) \cap S(\lambda^*|\bar L) \ne \{\textbf{0}_N\}$. Then, there exists a nonzero vector $t \in S(\lambda^*|L) \cap S(\lambda^*|\bar L)$. Let $\mu^*$ be an eigenvalue of $A+\lambda^*H$ with the corresponding eigenvector $\xi$. According to Lemma \ref{lemma1}, one has $\textbf{0}_{Nn} \ne t \otimes \xi \in \left\{{\mathop  \oplus \limits_{(i,j)\in  \Gamma(\mu^*)} {{  V}_{{i}{j}}}}\right\} \bigcap \left\{{\mathop  \oplus \limits_{(i,j)\in \bar \Gamma(\mu^*)} {{ \bar V}_{{i}{j}}}}\right\}$. Thus, it follows from Theorem \ref{theorem1} that the topological change $L \to \bar L$ is $\Psi$-indiscernible for the networked system (\ref{eq3})-(\ref{eq4}).
  \item For all $\lambda {\in} \sigma(L) {\cup} \sigma(\bar L)$, $S(\lambda|L) {\cap} S(\lambda|\bar L) {=} \{\textbf{0}_N\}$; But there exists $\lambda^* {\in} \sigma(L) \cup \sigma(\bar L)$ such that $\mathcal N(\Delta) \cap \{S(\lambda^*|L) \oplus S(\lambda^*|\bar L)\} \ne \{\textbf{0}_N\}$. It follows that there exists a nonzero vector $t\in \{S(\lambda^*|L) \oplus S(\lambda^*|\bar L)\}$ such that $\Delta t=\textbf{0}$. Let $t=t_1+t_2$, where $t_1 {\in} S(\lambda^*|L)$ and $t_2 {\in} S(\lambda^*|\bar L)$. Moreover,
   let $\mu^*$ be an eigenvalue of $A{+}\lambda^* H$, with the corresponding eigenvector $\xi$. It follows from Lemma \ref{lemma1} that $t_1{\otimes} \xi {\in} {\mathop  {\oplus} \limits_{(i,j){\in}  \Gamma(\mu^*)} {{  V}_{{i}{j}}}}$ and $t_2 {\otimes} \xi {\in} {\mathop  \oplus \limits_{(i,j)\in  \bar \Gamma(\mu^*)} {\bar {  V}_{{i}{j}}}}$. If $\left\{{\mathop  \oplus \limits_{(i,j)\in  \Gamma(\mu^*)} {{  V}_{{i}{j}}}}\right\} \bigcap \left\{{\mathop  \oplus \limits_{(i,j)\in \bar \Gamma(\mu^*)} {{ \bar V}_{{i}{j}}}}\right\} \ne \{\textbf{0}_{Nn}\}$, it follows from Theorem \ref{theorem1} that the topological change $L \to \bar L$ for the networked system (\ref{eq3})-(\ref{eq4}) is $\Psi$-indiscernible. If $\left\{{\mathop  \oplus \limits_{(i,j)\in  \Gamma(\mu^*)} {{  V}_{{i}{j}}}}\right\} \bigcap \left\{{\mathop  \oplus \limits_{(i,j)\in \bar \Gamma(\mu^*)} {{ \bar V}_{{i}{j}}}}\right\} = \{\textbf{0}_{Nn}\}$, then $t_1 \otimes \xi +t_2 \otimes \xi \in \left\{{\mathop  \oplus \limits_{(i,j)\in  \Gamma(\mu^*)} {{  V}_{{i}{j}}}}\right\} \bigoplus \left\{{\mathop  \oplus \limits_{(i,j)\in \bar \Gamma(\mu^*)} {{ \bar V}_{{i}{j}}}}\right\}$. Since $t {\ne} \textbf{0}$ and $\xi {\ne} \textbf{0}$, it follows that $t_1 {\otimes} \xi {+}t_2 {\otimes} \xi{=}t {\otimes} \xi {\ne} \textbf{0}$. Noting that $\Delta t{=}\textbf{0}$, one gets that $(\Delta {\otimes} C)(t_1 {\otimes} \xi{+}t_2 {\otimes} \xi)\\{=}(\Delta {\otimes} C)(t {\otimes} \xi){=}(\Delta t){\otimes}(C\xi){=}\textbf{0}$. Thus, $\textbf{0} {\ne} t_1 {\otimes} \xi {+}t_2 {\otimes} \xi \\{\in} \mathcal N (\Delta {\otimes} C) {\cap} \left\{\left\{{\mathop  {\oplus} \limits_{(i,j){\in}  \Gamma(\mu^*)} {{  V}_{{i}{j}}}}\right\}  {\bigoplus} \left\{{\mathop  {\oplus} \limits_{(i,j){\in} \bar \Gamma(\mu^*)} {{ \bar V}_{{i}{j}}}}\right\}\right\}$, which indicates that $\mathcal N (\Delta {\otimes} C) {\cap} \left\{\left\{{\mathop  {\oplus} \limits_{(i,j){\in}  \Gamma(\mu^*)} {{  V}_{{i}{j}}}}\right\}  {\bigoplus}\right.$ $\left. \left\{\mathop  {\oplus} \limits_{(i,j){\in} \bar \Gamma(\mu^*)} {{ \bar V}_{{i}{j}}}\right\}\right\} \ne \{\textbf{0}_{Nn}\}$. According to Theorem \ref{theorem1}, the topological change $L \to \bar L$ is $\Psi$-indiscernible for the networked system (\ref{eq3})-(\ref{eq4}).

\end{itemize}
Therefore, if the topological change $L \to \bar L$ is $\Delta$-indiscernible for system (\ref{eq11}), then this topological change is $\Psi$-indiscernible for the whole networked system (\ref{eq3})-(\ref{eq4}).
\hfill $\blacksquare$

The effectiveness of Corollary \ref{corollary6} is demonstrated by the
following example.

\begin{example}
Consider a simple network of three connected identical nodes, shown in (a) of Fig. \ref{fig1}, with $w_{21}=w_{32}=w_{23}=1$. Suppose that the output of the second node and the third node can be observed, i.e., $\delta_2=\delta_3=1$, with
\[A = \left[ {\begin{array}{*{20}{c}}
1&0\\
1&1
\end{array}} \right], \;H = \left[ {\begin{array}{*{20}{c}}
0&0\\
0&1
\end{array}} \right],\;C=\left[ {\begin{array}{*{20}{c}}
1&0\\
0&1
\end{array}} \right].\]

It can be easily verified that $L = \left[ {\begin{array}{*{20}{c}}
0&0&0\\
1&0&1\\
0&1&0
\end{array}} \right]$ and $\Delta=\left[ {\begin{array}{*{20}{c}}
0&1&0\\
0&0&1
\end{array}} \right]$. Assume that the edge from node $3$ to node $2$ is removed. Then, the new topology matrix is $\bar L = \left[ {\begin{array}{*{20}{c}}
0&0&0\\
1&0&0\\
0&1&0
\end{array}} \right]$.


It is easy to verify that $S(0|L)=span\{e_1-e_3\}$ and $S(0|\bar L)=span\{e_3\}$. It follows that $\mathcal N (\Delta) \cap \{S(0|L) \oplus S(0|\bar L)\} \ne \{\textbf{0}_3\}$, which implies that the topological change $L \to \bar L$ is $\Delta$-indiscernible for system (\ref{eq11}). Therefore, it follows from Corollary \ref{corollary6} that the topological change $L \to \bar L$ is $\Psi$-indiscernible for the whole networked system.

Actually, one cannot detect the topological change from the output trajectories when there exists some pair of initial states $(X_0,\;\bar X_0)$ satisfying $\Psi{e^{\Phi t}}{X_0} = \Psi{e^{\bar \Phi t}}{\bar X_0}$. An example is $X_0=[0,1,0,0,0,-1]^T$ and $\bar X_0=[0,0,0,0,0,1]^T$, with which the two networked systems will generate exactly the same output trajectories.

\end{example}

The above corollary reveals that the network topologies and the sensor locations have significant effect on the $\Psi$-discernibility for the whole networked system. In the following, the effect of the inner interactions, the node-system dynamics and output on the $\Psi$-discernibility is discussed.

\begin{corollary}
\label{corollary5}
Consider the networked system (\ref{eq3})-(\ref{eq4}). If a topological change $L \to \bar L$ is $\Psi$-discernible for the networked system (\ref{eq3})-(\ref{eq4}), then the following two conditions hold simultaneously:
\begin{enumerate}
  \item $(A,H)$ is observable;
  \item $(A+\lambda H, C)$ is observable, for all $\lambda \in \{\sigma(L) \cup \sigma(\bar L)\}$;
\end{enumerate}
\end{corollary}

\textbf{Proof:}
Using Theorem \ref{theorem1}, it can be proved easily. Thus, the proof is omitted.

\hfill $\blacksquare$

\begin{remark}
The possibility of detecting topological variations for networked LTI systems by observing network states has been investigated in \cite{haoyuqing}, with some lower-dimensional discernibility conditions established. However, complete observation of the full network states is unrealistic in practical applications. In most situations, only partial information about the node-system state is accessible, and only a subset of nodes are available for measurement. In this paper, the results in \cite{haoyuqing} are generalized to the case of detecting topological changes by observing output trajectories. Thus, the new conditions here are more general and have broader applicability in practice.
Compared with the conditions given in \cite{detecting_topology} and \cite{indiscernible}, which require the network topology to be undirected, the new conditions remove this requirement.
\end{remark}

\begin{remark}
The $\Psi$-discernibility of topological variations has been analyzed through eigenanalysis of the original and the modified networks in Section \ref{section3}. However, for most large-scale networked systems with higher-dimensional node-systems, this method cannot be applied efficiently. In this section, some lower-dimensional conditions on the $\Psi$-discernibility of topological changes are established. Only the properties of some {smaller matrices} are required to be checked. These lower-dimensional conditions allow to check the $\Psi$-discernibility {much more easily}, which also reveal how the topological variations, sensor locations, node-system dynamics and output, as well as inner interactions altogether affect the $\Psi$-discernibility of the topological variation.
\end{remark}

\section{{Output discernibility of topological variations for multi-agent systems}}
\label{section5}
In this section, the output discernibility of topological variations for multi-agent systems is revisited.
\subsection{Problem statement}
Consider a multi-agent system consisting of $N$ agents as follows:
 \[\left\{ \begin{array}{l}
{{\dot x}_i} = A{x_i} + B{u_i}, \\
{y_i} = C{x_i},
\end{array} \right.\;\;\;\;\;\;i = 1,2, \cdots ,N,\]
where $x_i \in \mathbb{R}^n$, $u_i \in \mathbb{R}^p$ and $y_i \in \mathbb{R}^m$ are the state, the input and the output of the $i$th agent, respectively; $A \in \mathbb{R}^{n\times n}$, $B \in \mathbb{R}^{n \times p}$ and $C \in \mathbb{R}^{m\times n}$ are the state matrix, the input matrix and the output matrix, respectively.

Agent $i$ is a neighbor of agent $j$ if its state is known by agent $j$. Here, assume that the neighboring relationships are fixed, which can be described by an undirected and weighted graph $\mathcal G=(\mathcal V, \mathcal E, \mathcal W)$. The coupling input $u_i$ to agent $i$ is determined by the diffusive coupling rule based on the neighboring relations as follows:

\[{u_i} = \sum\limits_{j \in {N_i}} {{w_{ij}}({x_j} - {x_i})}, \]
where $w_{ij}>0$ with $w_{ij}=w_{ji}$ for $(i,j) \in \mathcal E${,} and $N_i$ denotes the neighbor set of node $i$.

Let ${\mathcal{L}}=[l_{ij}] \in \mathbb{R}^{N \times N}$ be the graph Laplacian induced by $\mathcal G$, with
\[{l_{ij}}: = \left\{ \begin{array}{l}
\sum\limits_{k \in {N_i}} {{w_{ik}}} ,\;\;\;\;\;\;\;j = i;\\
 - {w_{ij}},\;\;\;\;\;\;\;\;\;\;j \ne i.
\end{array} \right.\]
Moreover, let $\mathcal Q \subset \mathcal V$ denote the subset of nodes
whose output is available for measurement, and $\Delta=col(e_i^T,\;i \in \mathcal Q) \in \mathbb{R}^{|\mathcal Q|\times N}$.
Let $X=[x_1^T, x_2^T,\cdots,x_N^T]^T$ and $Y=col(y_i,\;i \in \mathcal Q)$ be the state and the output of the whole multi-agent system, respectively. Then, the multi-agent system can be rewritten in a compact form as
\begin{equation}
\label{eq12}
\left\{ \begin{array}{l}
\dot X = FX,\\
Y=M X,
\end{array} \right.
\end{equation}
with
\begin{equation}
\label{eq13}
F=I_N \otimes A-{\mathcal{L}} \otimes B,\;
M=\Delta \otimes C.
\end{equation}

The nominal multi-agent system is represented by $({\mathcal G}, F, M)$. Here, consider a variation in the network structure, thereafter the new multi-agent system is described by
\begin{equation}
\label{eq14}
\left\{ \begin{array}{l}
\dot  {\bar X} = {\bar F} {\bar X},\\
\bar Y=M {\bar X},
\end{array} \right.
\end{equation}
with
\begin{equation}
\label{eq15}
\bar F=I_N \otimes A-{\bar{\mathcal{L}}} \otimes B,
\end{equation}
where ${\bar{\mathcal{L}}}$ is the graph Laplacian induced by the new topology $\bar{\mathcal G}$. Denote the new multi-agent system by $(\bar{\mathcal G}, \bar{F},M)$. Before moving on, the definitions of indiscernible initial state and always-discernible topological change are reviewed.

\begin{definition}\cite{indiscernible}
Consider the multi-agent system (\ref{eq12})-(\ref{eq13}). An initial state $X_0 \in \mathbb{R}^{Nn}$ is called indiscernible with respect
to the topological change $\mathcal L \to \bar {\mathcal L}$ if and only if
\[X(t) =e^{Ft}X_0=e^{\bar F t}X_0= \bar X(t),\;\;\forall t \ge 0.\]

\end{definition}

\begin{definition}\cite{indiscernible}
For the multi-agent system (\ref{eq12})-(\ref{eq13}), a topological change $\mathcal L \to \bar {\mathcal L}$
is called always-discernible if there is no (nontrivial) indiscernible
initial state. Otherwise, the topological change is called possibly-indiscernible.
\end{definition}

For the multi-agent systems $({\mathcal G}, F, M)$ and $(\bar{\mathcal G}, \bar{F},M)$, let the set of all the real $M$-indiscernible pairs of initial states be denoted as $T_M(F,\bar F)=\{(X_0,\bar X_0)\in \mathbb{R}^{Nn} \times \mathbb{R}^{Nn}| M e^{Ft}X_0=M e^{\bar Ft}\bar X_0, \forall t \ge 0\}$. Moreover, let $T_I(F,\bar F)=\{(X_0, X_0)\in \mathbb{R}^{Nn} \times \mathbb{R}^{Nn}| e^{Ft}X_0= e^{\bar Ft} X_0, \forall t \ge 0\}$. Noting that an indiscernible initial state always generates indiscernible output trajectories, one has $T_I(F,\bar F) \subseteq  T_M(F,\bar F)$. In what follows, the conditions under which the output matrix $M$ guarantees $T_M(F,\bar F)=T_I(F,\bar F)$ will be investigated. If $T_M(F,\bar F)=T_I(F,\bar F)$, the matrix $M$ is said to ensure the output discernibility of the multi-agent systems $({\mathcal G}, F, M)$ and $(\bar{\mathcal G}, \bar{F},M)$.

{\subsection{A counterexample and a {new} condition on output discernibility}}

Recall a condition on the output discernibility of the multi-agent systems $({\mathcal G}, F, M)$ and $(\bar{\mathcal G}, \bar{F},M)$ proposed in \cite{detecting_topology}, copied as follows{:}

{\begin{theorem}\label{theorem3}\cite{detecting_topology}
The matrix $M$ ensures the output discernibility of the multi-agent systems $(\mathcal G, F, M)$ and $( \bar {\mathcal G}, \bar F, M)$ if and only if the following conditions hold simultaneously:
\begin{enumerate}
  \item the pair $(\mathcal L,\Delta)$ is observable;
  \item the pair $(\bar {\mathcal L}, \Delta)$ is observable;
  \item the pair $(A-\lambda B,C)$ is observable for every $\lambda \in \sigma(\mathcal L) \cup \sigma(\bar {\mathcal L})$;
  \item the matrix $\Delta$ ensures the output discernibility of the two systems with state matrices $\mathcal L$ and $\bar {\mathcal L}$;
  \item the matrix $C$ ensures the output discernibility of the two systems with state matrices $A-\lambda B$ and $A-\bar \lambda B$, for $\lambda \in \sigma(\mathcal L)$ and $\bar \lambda \in \sigma(\bar {\mathcal L})$ with the corresponding eigenvectors $t$ and $\bar t$, respectively, which satisfy $\Delta t=\Delta \bar t$.
      \end{enumerate}
\end{theorem}}

In this theorem, it was claimed that the conditions are necessary and sufficient. However, the following counterexample shows that the conditions may not be sufficient.

\begin{example}
\label{example4}
Consider a multi-agent system consisting of three connected identical agents, shown in (a) of Fig. \ref{fig4}, with $w_{12}=w_{21}=w_{23}=w_{32}=w_{13}=w_{31}=1$. Suppose that the output of the first node and the third node can be observed, i.e., $\mathcal Q=\{1,\;3\}$, with
\[A = \left[ {\begin{array}{*{20}{c}}
1&1\\
0&2
\end{array}} \right],\;B = \left[ {\begin{array}{*{20}{c}}
0&0\\
0&1
\end{array}} \right],\;C = \left[ {\begin{array}{*{20}{c}}
1&0
\end{array}} \right].\]
\begin{figure}[htb]
\centering
\subfloat[${\mathcal G}$]{%
\includegraphics[width=.2\textwidth]{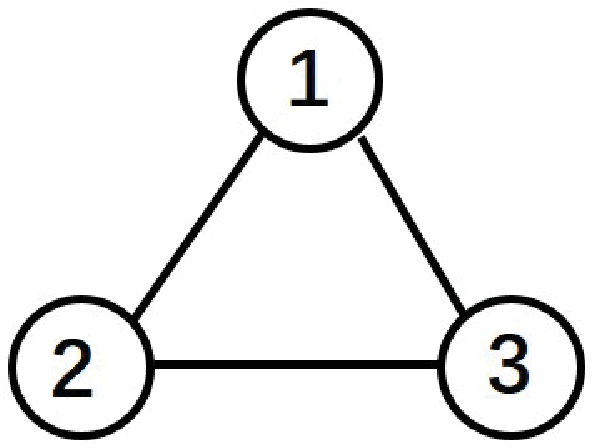}}
\subfloat[$\bar {\mathcal G}$]{%
\includegraphics[width=.25\textwidth]{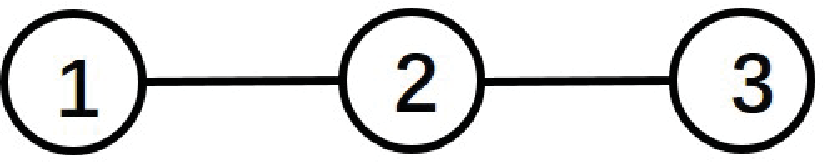}}\hfill
\caption{Network topologies}
\label{fig4}
\end{figure}

It can be easily verified that $\mathcal L = \left[ {\begin{array}{*{20}{c}}
2&{ - 1}&{ - 1}\\
{ - 1}&2&{ - 1}\\
{ - 1}&{ - 1}&2
\end{array}} \right]$ and $\Delta=\left[ {\begin{array}{*{20}{c}}
{ 1}&0&{ 0}\\
{ 0}&{ 0}&1
\end{array}} \right]$. Assume that the edge connecting nodes $3$ and $1$ is removed. Then, the new network is shown in (b) of Fig. \ref{fig4}, with the Laplacian matrix $\bar {\mathcal L} = \left[ {\begin{array}{*{20}{c}}
1&{ - 1}&0\\
{ - 1}&2&{ - 1}\\
0&{ - 1}&1
\end{array}} \right]$.

It is easy to verify that both $(\mathcal L,\Delta)$ and $(\bar {\mathcal L},\Delta)$ are observable. The eigenvalues of $\mathcal L$ are $\lambda_1=0$, $\lambda_2=1$ and $\lambda_3=3$. Moreover, the eigenvalues of $\bar {\mathcal L}$ are $\bar \lambda_1=0$ and $\bar \lambda_2=3$. One can easily verify that $(A-\lambda B,C)$ is observable, for $\lambda=0,\;1,\;3$.

Next, condition $(4)$ will be checked. The eigenspaces of ${\mathcal L}$ are $S(0|\mathcal L)=span\{e_1+e_2+e_3\}$, $S(1|\mathcal L)=span\{-e_1+e_3\}$ and $S(3|\mathcal L)=span\{e_1-2e_2+e_3\}$. Moreover, the eigenspaces of $\bar {\mathcal L}$ are $S(0|\bar {\mathcal L})=span\{e_1+e_2+e_3\}$ and $S(3|\bar {\mathcal L})=span\{-e_1+e_2,\;-e_1+e_3\}$. One can easily verify that
${T_\Delta(\mathcal L,\bar {\mathcal L})} = \left\{ (X_0,X_0)|X_0 \in span\left\{ e_1+e_2+e_3,\;e_1-2e_2+e_3 \right\}\right\}$. Moreover, it is obvious that $\mathcal L$ and $\bar {\mathcal L}$ have two common eigenpairs, which are $(0,e_1+e_2+e_3)$ and $(3,e_1-2e_2+e_3)$.
From the results given in \cite{indiscernible}, it can be verified that
$T_I(\mathcal L,\bar {\mathcal L}) = \left\{ (X_0,X_0)|X_0 \in span\left\{ e_1+e_2+e_3,\;e_1-2e_2\right.\right. $ $\left.\left.+e_3 \right\}\right\}=T_\Delta(\mathcal L,\bar {\mathcal L})$.
 Thus, matrix $\Delta$ ensures the output discernibility of the two systems with state matrices ${\mathcal L}$ and $\bar {\mathcal L}$.

Finally, condition $(5)$ will be checked. Since $t=e_1-2e_2+e_3 \in S(3|\mathcal L)$ and $\bar t=e_1+e_2+e_3 \in S(0|\bar {\mathcal L})$ satisfy that $\Delta t=\Delta \bar t$, one needs to check whether matrix $C$ can ensure the output discernibility of the two systems with state matrices $A-3B$ and $A$. The eigenspaces of $A-3B$ are $S(1|A-3B)=span\{e_1\}$ and $S(-1|A-3B)=span\{e_1-2e_2\}$. Moreover, the eigenspaces of $A$ are $S(1|A)=span\{e_1\}$ and $S(2|A)=span\{e_1+e_2\}$. It is easy to verify that $T_C(A-3B,A)=\left\{ (X_0,X_0)|X_0 \in span\{e_1\} \right\}=T_I(A-3B,A)$. Thus, matrix $C$ ensures the output discernibility of the two systems with state matrices $A-3B$ and $A$. Also, since $t=-e_1+e_3 \in S(1|\mathcal L)$ and $\bar t =-e_1+e_3 \in S(3|\bar {\mathcal L})$ satisfy that $\Delta t=\Delta \bar t$, one needs to check whether matrix $C$ can ensure the output discernibility of the two systems with state matrices $A-B$ and $A-3B$. Note that the eigenspace of $A-B$ is $S(1|A-B)=span\{e_1\}$. It is easy to verify that $T_C(A-B,A-3B)=\left\{ (X_0,X_0)|X_0 \in span\{e_1\} \right\}=T_I(A-B,A-3B)$. Thus, matrix $C$ ensures the output discernibility of the two systems with state matrices $A-B$ and $A-3B$. Therefore, condition $(5)$ holds.

From the results given in \cite{detecting_topology}, it would be deduced that $M$ ensures the output discernibility of the multi-agent systems $(\mathcal G, F, M)$ and $( \bar {\mathcal G}, \bar F, M)$.
However, if one chooses $(X_0,\;\bar X_0)=([0,0,3,0,0,0]^T,\;\textbf{0}_6) \notin T_I(F,\bar F)$, the two multi-agent systems  generate exactly the same output trajectories. Thus, $T_M(F,\bar F) \ne T_I(F,\bar F)$, which indicates that $M$ actually does not ensure the output discernibility.
Therefore, the sufficiency of the condition given in \cite{detecting_topology} does not hold.

\end{example}

%
%

Note that Theorem \ref{theorem10} can also be used to verify the $M$-discernibility of the topological variation ${\mathcal L} \to \bar {\mathcal L}$ for the multi-agent system (\ref{eq12})-(\ref{eq13}). Moreover,
a topological variation is $M$-discernible if and only if the topological variation is always-discernible and matrix $M$ ensures the output discernibility. The first condition in Theorem \ref{theorem10} requires that the topological change ${\mathcal L} \to \bar {\mathcal L}$ is always-discernible, while the second one requires that matrix $M$ ensures the output discernibility. Based on the results in Section \ref{section3}, a new necessary and sufficient condition on the output discernibility of the multi-agent systems $({\mathcal G}, F, M)$ and $(\bar{\mathcal G}, \bar{F},M)$ can be established as follows.

Let $\phi$ be an eigenvalue of $F$, with the corresponding eigenspace
$S(\phi|F)$.
Moreover, let $\bar \phi$ be an eigenvalue of $\bar F$, with the associated eigenspace
$S(\bar \phi|\bar F)$. Now, a condition on the output discernibility of the multi-agent systems $({\mathcal G}, F, M)$ and $(\bar{\mathcal G}, \bar{F},M)$ is presented.

\begin{corollary}\label{corollary8}
The matrix $M$ ensures the output discernibility of the multi-agent systems $({\mathcal G}, F, M)$ and $(\bar{\mathcal G}, \bar{F},M)$ if and only if for all $\phi \in \sigma(F)\cup \sigma(\bar F)$, one has
\[\mathcal N (M) \cap \{S(\phi|F) + S(\phi|\bar F)\}=\{\textbf{0}_{Nn}\}.\]
\end{corollary}

\begin{remark}
Some conditions on the output discernibility of topological variations for undirected networks have been established in \cite{detecting_topology}. The new condition proposed in Corollary \ref{corollary8} can handle both undirected and directed networks, thus is more general.

\end{remark}

By using Lemma \ref{lemma1}, the eigenspaces of $F$ and $\bar  F$ can be expressed through the generalized eigenvectors of some matrices with lower dimensions. Similarly as in Section \ref{section4}, some lower-dimensional conditions on the output discernibility of the multi-agent systems $({\mathcal G}, F, M)$ and $(\bar{\mathcal G}, \bar{F},M)$ can also be established.

Finally, the effectiveness of the new condition can also be illustrated by Example \ref{example4}.

It suffices to observe that $S(1|F)=span\{(e_1+e_2+e_3)\otimes e_1,\;(-e_1+e_3)\otimes e_1,\;(e_1-2e_2+e_3) \otimes e_1\}$. There exists $\eta=[0,0,3,0,0,0]^T =[0,0,3,0,0,0]^T+\textbf{0}_6 \in S(1|F)+S(1|\bar F)$ such that $M\eta=0$. Thus, $\mathcal N (M) \cap \{S(1|F) + S(1|\bar F)\}\ne\{\textbf{0}_{Nn}\}$. From Corollary \ref{corollary8}, it follows that the matrix $M$ does not ensure the output discernibility of the multi-agent systems $({\mathcal G}, F, M)$ and $(\bar{\mathcal G}, \bar{F},M)$.

\section{Conclusions}
\label{section6}
This paper has investigated the conditions under which a topological variation in networked LTI systems can be detected by observing output trajectories. The considered network topology can be general, directed and weighted. A necessary and sufficient condition on the $\Psi$-discernibility of topological variations has been established in terms of the eigenspaces of the original and the modified networks. Further, by taking the network structures into account, some lower-dimensional conditions on the $\Psi$-discernibility have been derived. The conditions have generalized the results given in \cite{haoyuqing}, assuming only partial state variables of the networks are available for measurement. Moreover, the output discernibility of topological variations for multi-agent systems has been revisited. It is found that the sufficiency of the criterion given in \cite{detecting_topology} does not hold. Consequently, a complete necessary and sufficient condition is established.
In future studies, the important issue of restoring the discernibility of topological variations for networked LTI systems will be investigated.

\end{document}